\documentclass[11pt]{amsart}
\usepackage{euscript}
\usepackage{amssymb}
\usepackage{amsmath}
\usepackage{epic}

\numberwithin{equation}{section}

\newtheoremstyle{my}{1.5em}{0.5em}{\em}{}{\sc}{.}{0.5em}{}

\theoremstyle{my}

\theoremstyle{my}
\newtheorem{thm}{Theorem}[section]
\newtheorem{Theorem}[thm]{Theorem}
\newtheorem*{Theorem*}{Theorem}
\newtheorem{Corollary}[thm]{Corollary}

\newtheorem*{corollary*}{Corollary}
\newtheorem{Lemma}[thm]{Lemma}

\newtheorem{Conjecture}[thm]{Conjecture}
\newtheorem*{conjecture*}{Conjecture}

\newtheorem*{question*}{Question}

\newtheorem*{definitions*}{Definitions}

\newtheorem*{rem*}{Remark}
\newtheorem{Remark}[thm]{Remark}

\newtheorem*{remark*}{Remark}

\newtheorem*{remarks*}{Remarks}
\newtheorem*{example*}{Example}
\newtheorem{Example}[thm]{Example}
\newtheorem*{examples*}{Examples}

\newtheorem*{convention*}{Convention}
\newtheorem*{conventions*}{Conventions}
\newtheorem*{Note*}{Note}
\newtheorem*{exercise*}{Exercise}
\newtheorem*{bibliographical-note*}{Bibliographical note}

\newcommand{\Acknowledgements}{{\em Acknowledgements.} }

\newcommand{\R}{\mathbb{R}}
\newcommand{\Z}{\mathbb{Z}}
\newcommand{\Q}{\mathbb{Q}}
\newcommand{\C}{\mathbb{C}}

\renewcommand{\P}{\mathbb{P}}
\newcommand{\half}{{\textstyle\frac{1}{2}}}

\newcommand{\iso}{\cong}           
\newcommand{\smooth}{C^\infty}

\newcommand{\im}{\mathrm{im}}

\newcommand{\End}{\mathrm{End}}

\newcommand{\K}{\mathbb{K}}
\newcommand{\A}{\mathcal{A}}
\newcommand{\F}{\mathcal{F}}
\newcommand{\W}{\mathcal{W}}

\title[Cotangent bundles]{The symplectic geometry of cotangent bundles from a categorical viewpoint}
\author{Kenji Fukaya, \, Paul Seidel, \, Ivan Smith}

\begin{document}
\maketitle

\begin{abstract}
We describe various approaches to understanding Fukaya categories of cotangent bundles.  All of the approaches rely on introducing a suitable class of noncompact Lagrangian submanifolds. We review the work of Nadler-Zaslow \cite{NZ,Nadler} and the authors \cite{FSS}, before discussing a new approach using family Floer cohomology \cite{Fukaya} and the ``wrapped Fukaya category''. The latter, inspired by Viterbo's symplectic homology, emphasises the connection to loop spaces, hence seems particularly suitable when trying to extend the existing theory beyond the simply-connected case.
\end{abstract}

A classical problem in symplectic topology is to describe exact Lagrangian submanifolds inside a cotangent bundle. The main conjecture, usually attributed to Arnol'd \cite{Arnold}, asserts that any (compact) submanifold of this kind should be Hamiltonian isotopic to the zero-section. In this sharp form, the result is known only for $S^2$ and $\R\P^2$, and the proof uses methods which are specifically four-dimensional (both cases are due to Hind \cite{Hind}; concerning the state of the art for surfaces of genus $>0$, see \cite{Hind-Ivrii}). In higher dimensions, work has concentrated on trying to establish topological restrictions on exact Lagrangian submanifolds. There are many results dealing with assorted partial aspects of this question
(\cite{LS,Vit:loop,Vit,Buhovsky,Seidel} and others), using a variety of techniques. Quite recently, more categorical methods have been added to the toolkit, and these have led to a result covering a fairly general situation.
The basic statement (which we will later generalise somewhat) is as follows:

\begin{Theorem}[Nadler, Fukaya-Seidel-Smith] \label{Thm:main}
Let $Z$ be a closed, simply connected manifold which is spin, and $M = T^*Z$ its cotangent bundle. Suppose that $L \subset M$ is an exact closed Lagrangian submanifold, which is also spin, and additionally has vanishing Maslov class. Then: (i) the projection $L \hookrightarrow M \rightarrow Z$ has degree $\pm 1$; (ii) pullback by the projection is an isomorphism $H^*(Z;\K) \iso H^*(L;\K)$ for any coefficient field $\K$; and (iii) given any two Lagrangian submanifolds $L_0,L_1 \subset M$ with these properties, which meet transversally, one has $|L_0 \cap L_1| \geq \mathrm{dim}\, H^*(Z;\K)$.
\end{Theorem}

Remarkably, three ways of arriving at this goal have emerged, which are essentially independent of each other, but share a basic philosophical outlook. One proof is due to Nadler \cite{Nadler}, building on earlier work of Nadler and Zaslow \cite{NZ} (the result in \cite{Nadler} is formulated for $\K = \C$, but it seems that the proof goes through for any $\K$). Another one is given in \cite{FSS}, and involves, among other things, tools from \cite{FCPLT} (for technical reasons, this actually works only for $\mathrm{char}(\K) \neq 2$). The third one, which is collaborative work of the three authors of this paper, is not complete at the time of writing, mostly because it relies on ongoing developments in general Floer homology theory. In spite of this, we included a description of it, to round off the overall picture.

The best starting point may actually be the end of the proof, which can be taken to be roughly the same in all three cases. Let $L$ be as in Theorem \ref{Thm:main}. The Floer cohomology groups $HF^*(T^*_x,L)$, where $T^*_x \subset M$ is the cotangent fibre at some point $x \in Z$, form a flat bundle of $\Z$-graded vector spaces over $Z$, which we denote by $E_L$. There is a spectral sequence converging to $HF^*(L,L) \iso H^*(L;\K)$, whose $E_2$ page is
\begin{equation} \label{eq:ss}
E_2^{rs} = H^r(Z;End^s(E_L)),
\end{equation}
$End^*(E_L) = Hom^*(E_L,E_L)$ being the graded endomorphism vector
bundle. Because of the assumption of simple connectivity of $Z$, $E_L$ is actually trivial, so the $E_2$ page is a ``box'' $H^*(Z) \otimes End(HF^*(T^*_x,L))$. The $E_2$ level differential goes from $(r,s)$ to $(r+2,s-1)$:
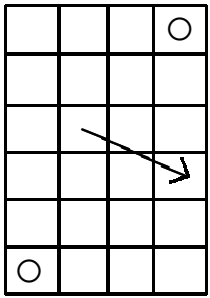
\begin{center}
\setlength{\unitlength}{0.00041667in}
\begingroup\makeatletter\ifx\SetFigFont\undefined%
\gdef\SetFigFont#1#2#3#4#5{%
  \reset@font\fontsize{#1}{#2pt}%
  \fontfamily{#3}\fontseries{#4}\fontshape{#5}%
  \selectfont}%
\fi\endgroup%
{\renewcommand{\dashlinestretch}{30}
\begin{picture}(2594,3734)(0,-10)
\thicklines
\drawline(1322,3697)(1322,22)
\put(322,322){\circle{272}}
\drawline(997,2122)(2347,1522)
\drawline(997,2122)(2347,1522)
\drawline(2272,1747)(2347,1522)(2122,1447)
\drawline(2272,1747)(2347,1522)(2122,1447)
\drawline(22,1227)(2572,1227)
\drawline(22,1822)(2572,1822)
\drawline(22,2422)(2572,2422)
\drawline(22,3077)(2572,3077)
\drawline(697,3697)(697,22)
\drawline(22,622)(2572,622)
\drawline(2572,3697)(22,3697)(22,22)
	(2572,22)(2572,3697)
\drawline(1897,3697)(1897,22)
\put(2234,3397){\circle{272}}
\end{picture}
}
\end{center}
and similarly for the higher pages. Hence, the bottom left and top right corners of the box necessarily survive to $E_\infty$. Just by looking at their degrees, it follows that $HF^*(T^*_x,L) \iso \K$ must be one-dimensional (and, we may assume after changing the grading of $L$, concentrated in degree $0$). Given that, the spectral sequence degenerates, yielding $H^*(L;\K) \iso H^*(Z;\K)$. On the other hand, we have also shown that the projection $L \rightarrow Z$ has degree $\pm 1 = \pm \chi(HF^*(T^*_x,L))$. This means that the induced map on cohomology is injective, hence necessarily an isomorphism. Finally, there is a similar spectral sequence for a pair of Lagrangian submanifolds $(L_0,L_1)$, which can be used to derive the last part of Theorem \ref{Thm:main}.

At this point, we already need to insert a cautionary note. Namely, the approach in \cite{FSS} leads to a spectral sequence which only approximates the one in \eqref{eq:ss} (the $E_1$ term is an analogue of the expression above, replacing the cohomology of $Z$ by its Morse cochain complex, and the differential is only partially known). In spite of this handicap, a slightly modified version of the previous argument can be carried out successfully. The other two strategies (\cite{Nadler} and the unpublished approach) do not suffer from this deficiency, since they directly produce \eqref{eq:ss} in the form stated above.

From the description we have just given, one can already infer one
basic philosophical point, namely the interpretation of Lagrangian
submanifolds in $M$ as (some kind of) sheaves on the base $Z$. This
can be viewed as a limit of standard ideas about Lagrangian torus
fibrations in mirror symmetry \cite{FukayaAV,Kontsevich-Soibelman},
where the volume of the tori becomes infinite (there is no
algebro-geometric mirror of $M$ in the usual sense, so we borrow only
half of the mirror symmetry argument). The main problem is to prove
that the sheaf-theoretic objects accurately reflect the Floer
cohomology groups of Lagrangian submanifolds, hence in particular
reproduce $HF^*(L,L) \iso H^*(L;\K)$.  Informally speaking, this is
ensured by providing a suitable ``resolution of the diagonal'' in the
Fukaya category of $M$, which reduces the question to one about
cotangent fibres $L = T^*_x$. In saying that, we have implicitly
already introduced an enlargement of the ordinary Fukaya category,
namely one which allows noncompact Lagrangian submanifolds. There are several
possible ways of treating such submanifolds, leading to categories
with substantially different properties. This is where the three
approaches diverge:

\textbf{Characteristic cycles.} \cite{NZ} considers a class of
Lagrangian submanifolds which, at infinity, are invariant under
rescaling of the cotangent fibres (or more generally, asymptotically
invariant). Intersections at infinity are dealt with by small
perturbations (in a distinguished direction given by the normalized
geodesic flow; this requires the choice of a real analytic structure on $Z$).
An important source of inspiration is Kashiwara's construction
\cite{Kashiwara} of 
characteristic cycles for constructible sheaves on $Z$; and indeed,
Nadler proves 
that, once derived, the resulting version of the Fukaya category is
equivalent to the constructible derived category.  (A similar point of
view was taken in earlier papers of Kasturirangan and Oh
\cite{Kasturirangan-Oh,Kasturirangan-Oh2, Oh}).  Generally
speaking, to get a finite resolution of the diagonal, this category
has to be modified further, by restricting the behaviour at infinity;
however, if one is only interested in applications to closed
Lagrangian submanifolds, this step can be greatly simplified.

\textbf{Lefschetz thimbles.} The idea in \cite{FSS} is to embed the
Fukaya category of $M$ into the Fukaya category of a Lefschetz
fibration $\pi: X \rightarrow \C$. The latter class of categories is
known to admit full exceptional collections, given by any basis of
Lefschetz thimbles. Results from homological algebra (more precisely, the theory of mutations, see for
instance \cite{Goro}) then ensure the existence of a resolution of the
diagonal, in terms of Koszul dual bases. To apply this machinery one
has to construct a Lefschetz fibration, with an antiholomorphic
involution, whose real part $\pi_\R$ is a Morse function on $X_\R =
Z$. This can be done easily, although not in a canonical way, by using
techniques from real algebraic geometry. Roughly speaking, the
resulting Fukaya category looks similar to the category of sheaves
constructible with respect to the stratification given by the unstable
manifolds of $\pi_\R$, compare \cite{Kapranov}. However, because the
construction of $X$ is not precisely controlled, one does not expect
these two categories to agree. Whilst this is not a problem for the
proof of Theorem \ref{Thm:main}, it may be {\ae}sthetically
unsatisfactory. One possibility for improving the situation would be
to find a way of directly producing a Lefschetz fibration on the
cotangent bundle; steps in that direction are taken in \cite{Johns}.

\textbf{Wrapping at infinity.} The third approach remains within $M$,
and again uses Lagrangian submanifolds which are scaling-invariant at
infinity. However, intersections at infinity are dealt with by flowing
along the (not normalized) geodesic flow, which is a large
perturbation. For instance, after this perturbation, the intersections
of any two fibres will be given by all geodesics in $Z$ connecting the
relevant two points. In contrast to the previous constructions, this one is intrinsic to the differentiable manifold $Z$, and does not require a real-analytic or real-algebraic structure (there are of course technical choices to be made, such as the Riemannian metric and other perturbations belonging to standard pseudo-holomorphic curve theory; but the outcome is independent of those up to quasi-isomorphism).

Conjecturally, the resulting ``wrapped Fukaya category'' is equivalent to a full subcategory of the category of modules over $C_{-*}(\Omega Z)$, the dg (differential graded) algebra
of chains on the based loop space (actually, the Moore loop space,
with the Pontryagin product). Note that the classical bar construction
establishes a relation between this and the dg algebra of cochains on
$Z$; for
$\K = \R$, one can take this to be the algebra of differential forms; for
$\K = \Q$, it could be Sullivan's model; and for general $\K$ one
can use singular or Cech cohomology. If $Z$ is simply
connected, this relation leads to an equivalence of suitably defined
module categories, and one can recover \eqref{eq:ss} in this way. In fact, we propose a more geometric version of this argument, which involves an explicit functor from the wrapped Fukaya category to the category of modules over a dg algebra of Cech cochains.

\begin{Remark}
  \label{Rem:1}
  It is interesting to compare (\ref{eq:ss}) with the result of a
  naive geometric argument.  Suppose $L$ is closed and exact; under
fibrewise scaling $L \mapsto cL$, as $c\rightarrow \infty$,  $cL$
converges (in compact subsets) to a disjoint union of cotangent fibres
$\bigcup_{x\in L\cap i(Z)}
\,T_x^*$, where $i:Z\hookrightarrow M$ denotes the zero-section.  In
particular, for  large $c$ there is a canonical bijection
between points of $L\cap cL$ and points of $\{L\cap T_x^* \, | \,
x\in L\cap i(Z)\}$. Starting from this identification, and filtering by
energy, one expects to obtain a spectral sequence
\begin{equation}
\bigoplus_{x \in L\cap i(Z)} HF(L,T_x^*) \Rightarrow H^*(L)
\end{equation}
using exactness to identify
$HF(L,cL)\cong H^*(L)$.  This would (re)prove that $L\cap
i(Z)\neq \emptyset$ and that $HF(L,T_x^*)\neq 0$; and, in an informal
fashion, this provides  motivation for believing that
$L$ is ``generated by the fibres''. It seems hard,
however, to control the homology class of $L$ starting from this, because
it seems hard to gain sufficient control over $L\cap
i(Z)$.
\end{Remark}

The following three sections of the paper are each devoted to explaining one of these approaches. Then, in a concluding section, we take a look at the non-simply-connected case. First of all, there is a useful trick involving the spectral sequence \eqref{eq:ss} and finite covers of the base $Z$. In principle, this trick can be applied to any of the three approaches outlined above, but at the present state of the literature, the necessary prerequisites have been fully established only for the theory from \cite{NZ}. Applying that, one arrives at the following consequence, which appears to be new:

\begin{Corollary} \label{Cor:main}
The assumption of simple-connectivity of $Z$  can be removed from Theorem \ref{Thm:main}. This means that for all closed spin manifolds $Z$, and all exact Lagrangian submanifolds $L \subset T^*Z$ which are spin and have zero Maslov class, the conclusions (i)--(iii) hold.
\end{Corollary}

From a more fundamental perspective, the approach via wrapped Fukaya
categories seems particularly suitable for investigating cotangent
bundles of non-simply-connected manifolds, since it retains
information that is lost when passing from chains on $\Omega Z$ to
cochains on $Z$. We end by describing what this would mean (modulo one
of our conjectures, \ref{Conj:generate}) in the special case when $Z$
is a $K(\Gamma,1)$.   In the special case of the torus $Z=T^n$, Conjecture \ref{Conj:generate} can be sidestepped by direct geometric arguments, at least when $\mathrm{char}(\K)\neq 2$.  Imposing that condition, one finds that an arbitrary exact, oriented and spin Lagrangian submanifold $L\subset T^*T^n = (\C^*)^n$ satisfies the conclusions of Theorem \ref{Thm:main}, with no assumption on the Maslov class.

Finally, it is worth pointing out that for any oriented and spin Lagrangian submanifold $L\subset T^*Z$, and any closed spin manifold $Z$, the theory produces a $\Z/2\Z$-graded spectral sequence
\begin{equation}\label{eq:ungraded}
H(Z; \End(E_L)) \Rightarrow H(L)
\end{equation}
which has applications in its own right, for instance to the classification of ``local Lagrangian knots" in Euclidean space.  Eliashberg and Polterovich \cite{EP} proved that any exact Lagrangian $L\subset \C^2$ which co-incides with the standard $\R^2$ outside a compact subset is in fact Lagrangian isotopic to $\R^2$.  Again, their proof relies on exclusively four-dimensional machinery.  

\begin{Corollary} \label{Cor:euclidean}
Let $L\subset \C^n$ be an exact Lagrangian submanifold which co-incides with $\R^n$ outside a compact set.  Suppose that $L$ is oriented and spin.  Then (i) $L$ is acyclic and (ii) $\pi_1(L)$ has no non-trivial finite-dimensional complex representations.
\end{Corollary}

Using a result of Viterbo \cite{Vit} and standard facts from 3-manifold topology, when $n=3$ this implies that an oriented exact $L\subset \C^3$ which co-incides with $\R^3$ outside a compact set is actually diffeomorphic to $\R^3$.  To prove Corollary \ref{Cor:euclidean}, one embeds the given $L$ (viewed in a Darboux ball) into the zero-section of $T^*S^n$, obtaining an exact Lagrangian submanifold $L'$ of the latter which co-incides with the zero-section on an open set.  This last fact immediately implies that $E_{L'}$ has rank one, so we see $\End(E_{L'})$ is the trivial $\K$-line bundle, and the sequence (\ref{eq:ungraded}) implies $\textrm{rk}_{\K} H^*(S^n)\geq \textrm{rk}_{\K} H^*(L')$. On the other hand, projection $L' \rightarrow S^n$ is (obviously) degree $\pm 1$, which gives the reverse inequality.  Going back to $L\subset L'$, we deduce that $L$ is acyclic with $\K$-coefficients.  Since $\K$ is arbitrary, one deduces the first part of the Corollary. At this point, one knows \emph{a posteori} that the Maslov class of $L$ vanishes; hence so does that of $L'$, and the final statement of the Corollary then follows from an older result of Seidel \cite{Seidel}. When $n=3$, it follows that $L'$ is simply connected or a $K(\pi,1)$, at which point one can appeal to Viterbo's work.  (It is straightforward to deduce the acyclicity over fields of characteristic other than two in the other frameworks, for instance that of \cite{FSS}, but then the conclusion on the Maslov class does not follow.)

\Acknowledgements I.S. is grateful to Kevin Costello for encouragement and many helpful conversations. Some of this material was presented at the ``Dusafest'' (Stony Brook, Oct 2006). P.S. was partially supported by NSF grant DMS-0405516. I.S. acknowledges EPSRC grant EP/C535995/1.

\section{Constructible sheaves\label{sec:nz}}

This section should be considered as an introduction to the two papers \cite{NZ,Nadler}. Our aim is to present ideas from those papers in a way which is familiar to symplectic geometers. With that in mind, we have taken some liberties in the presentation, in particular omitting the (nontrivial) technical work involved in smoothing out characteristic cycles.

\subsection{Fukaya categories of Weinstein manifolds} \label{subsec:epsilon}
Let $M$ be a Weinstein manifold which is of finite type and
complete. Recall that a symplectic manifold $(M,\omega)$ is Weinstein
if it comes with a distinguished Liouville (symplectically expanding)
vector field $Y$, and a proper bounded below function $h: M
\rightarrow \R$, such that $dh(Y)$ is positive on a sequence of level
sets $h^{-1}(c_k)$, with $\lim_k c_k = \infty$. The stronger finite
type assumption is that $dh(Y)>0$ outside a compact subset of
$M$. Finally, completeness means that the flow of $Y$ is defined for
all times (for negative times, this is automatically true, but for
positive times it is an additional constraint). Note that the Liouville vector field defines a one-form $\theta = i_Y\omega$ with $d\theta = \omega$. At infinity, $(M,\theta)$ has the form $([0;\infty) \times N, e^r\alpha)$, where $N$ is a contact manifold with contact one-form $\alpha$, and $r$ is the radial coordinate. In other words, the end of $M$ is modelled on the positive half of the symplectization of $(N,\alpha)$. The obvious examples are cotangent bundles of closed manifolds, $M = T^*Z$, where $Y$ is the radial rescaling vector field, and $N$ the unit cotangent bundle.

We will consider exact Lagrangian submanifolds $L \subset M$ which are Legendrian at infinity. By definition, this means that $\theta|L$ is the derivative of some compactly supported function on $L$. Outside a compact subset, any such $L$ will be of the form $[0;\infty) \times K$, where $K \subset N$ is a Legendrian submanifold. Now let $(L_0,L_1)$ be two such submanifolds, whose structure at infinity is modelled on $(K_0,K_1)$. To define their Floer cohomology, one needs a way of resolving the intersections at infinity by a suitable small perturbation. The details may vary, depending on what kind of Legendrian submanifolds one wants to consider. Here, we make the assumption that $(N,\alpha)$ is real-analytic, and allow only those $K$ which are real-analytic submanifolds. Then,

\begin{Lemma} \label{Th:daft}
Let $(\phi^t_R)$ be the Reeb flow on $N$. For any pair $(K_0,K_1)$, there is an $\epsilon>0$ such that $\phi^t_R(K_0) \cap K_1 = \emptyset$ for all $t \in (0,\epsilon)$.
\end{Lemma}

This is a consequence of the Curve Selection Lemma \cite[Lemma
3.1]{Milnor}, compare \cite[Lemma 5.2.5]{NZ}. Recall that, when
defining the Floer cohomology of two Lagrangian submanifolds, one
often adds a Hamiltonian perturbation $H \in \smooth(M,\R)$ (for
technical reasons, this Hamiltonian is usually also taken to be
time-dependent, but we suppress that here). The associated Floer
cochain complex is generated by the flow lines $x: [0;1] \rightarrow
M$ of $H$ going from $L_0$ to $L_1$; equivalently, these are the
intersection points of $\phi^1_X(L_0) \cap L_1$, where $X$ is the
Hamiltonian vector field of $H$. We denote this cochain complex by
\begin{equation}
CF^*(L_0,L_1;H) = CF^*(\phi_X^1(L_0),L_1).
\end{equation}
In our case, we take an $H$ which at infinity is of the form $H(r,y) = h(e^r)$, where $h$ is a function with $h' \in (0;\epsilon)$. Then, $X$ is $h'(e^r)$ times the Reeb vector field $R$, hence $\phi^1_X(L_0) \cap L_1$ is compact by Lemma \ref{Th:daft}. Standard arguments show that the resulting Floer cohomology group $HF^*(L_0,L_1) = HF^*(L_0,L_1;H)$ is independent of $H$.  It is also invariant under compactly supported (exact Lagrangian) isotopies of either $L_0$ or $L_1$. Note that in the case where $K_0 \cap K_1 = \emptyset$, one can actually set $h = 0$, which yields Floer cohomology in the ordinary (unperturbed) sense. Finally, for $L_0 = L_1 = L$ one has the usual reduction to Morse theory, so that $HF^*(L,L) \iso H^*(L)$, even for noncompact $L$.

At this point, we need to make a few more technical remarks. For simplicity, all our Floer cohomology groups are with coefficients in some field ${\mathbb K}$. If $char({\mathbb K}) \neq 2$, one needs (relative) $Spin$ structures on all Lagrangian submanifolds involved, in order to address the usual orientation problems for moduli spaces \cite{FO3}. Next, Floer cohomology groups are generally only $\Z/2$-graded. One can upgrade this to a $\Z$-grading by requiring that $c_1(M) = 0$, and choosing gradings of each Lagrangian submanifold. For the moment, we do not need this, but it becomes important whenever one wants to make the connection with objects of classical homological algebra, as in Theorem \ref{th:nz1} below, or in \eqref{eq:ss}.

\begin{Example} \label{Ex:cotangent}
Consider the case of cotangent bundles $M = T^*Z$ (to satisfy the general requirements above, we should impose real-analyticity conditions, but that is not actually necessary for the specific computations we are about to do). A typical example of a Lagrangian submanifold $L \subset M$ satisfying the conditions set out above is the conormal bundle $L = \nu^*W$ of a closed submanifold $W \subset Z$. If $(L_0,L_1)$ are conormal bundles of transversally intersecting submanifolds $(W_0,W_1)$, then
\begin{equation}
HF^*(L_0,L_1) \iso H^{*-codim(W_0)}(W_0 \cap W_1).
\end{equation}
This is easy to see (except perhaps for the grading), since the only intersection points of the $L_k$ lie in the zero-section. All of them have the same value of the action functional, and standard Morse-Bott techniques apply.

As a parallel but slightly different example, let $W \subset Z$ be an
open subset with smooth boundary. Take a function $f: \overline{W}
\rightarrow \R$ which is strictly positive in the interior, zero on
the boundary, and has negative normal derivative at all boundary
points. We can then consider the graph of $d(1/f)$, which is a
Lagrangian submanifold of $M$, asymptotic to the positive part of the
conormal bundle of $\partial W$. By a suitable isotopy, one can deform
the graph so that it agrees at infinity with that conormal
bundle. Denote the result by $L$. Given two such subsets $W_k$ whose
boundaries intersect transversally, one then has
\cite{Kasturirangan-Oh, Kasturirangan-Oh2, Oh}
\begin{equation}
HF^*(L_0,L_1) \iso H^*(W_1 \cap \overline{W}_0,W_1 \cap \partial W_0).
\end{equation}
Note that in both these cases, the Lagrangian submanifolds under consideration do admit natural gradings, so the isomorphisms are ones of $\Z$-graded groups.
\end{Example}

We will need multiplicative structures on $HF^*$, realized on the chain level by an $A_\infty$-category structure. The technical obstacle, in the first nontrivial case, is that the natural triangle product
\begin{equation} \label{eq:strict-product}
\begin{aligned}
&
CF^*(L_1,L_2;H_{12}) \otimes CF^*(L_0,L_1;H_{01}) \\ & = CF^*(\phi_{X_{12}}^1(L_1),L_2) \otimes
CF^*(\phi_{X_{12}}^1\phi_{X_{01}}^1(L_0),\phi_{X_{12}}^1(L_1)) \\ & \longrightarrow CF^*(\phi_{X_{12}}^1
\phi_{X_{01}}^1(L_0),L_2)
\end{aligned}
\end{equation}
does not quite land in $CF^*(L_0,L_2;H_{02}) = CF^*(\phi_{X_{02}}^1(L_0),L_2)$. For instance, if one takes the same $H$ for all pairs of Lagrangian submanifolds, the output of the product has $\phi_X^2 = \phi_{2X}^1$ instead of the desired $\phi_X^1$. The solution adopted in \cite{NZ} is (roughly speaking) to choose all functions $h$ involved to be very small, in which case the deformation from $X$ to $2X$ induces an actual isomorphism of Floer cochain groups. The downside is that this can only be done for a finite number of Lagrangian submanifolds, and more importantly, for a finite number of $A_\infty$-products at a time. Hence, what one gets is a partially defined $A_d$-structure (for $d$ arbitrarily large), from which one then has to produce a proper $A_\infty$-structure; for some relevant algebraic results, see \cite[Lemma 30.163]{FO3}. As an alternative, one can take all functions $h$ to satisfy $h(t) = \log\,t$, which means $H(r,y) = r$. Then, $\phi_X^2$ is conjugate to a compactly supported perturbation of $\phi_X^1$ (the conjugating diffeomorphism is the Liouville flow $\phi_Y^t$ for time $t = -\log(2)$). By making the other choices in a careful way, one can then arrange that \eqref{eq:strict-product} takes values in a Floer cochain group which is isomorphic to $CF^*(L_0,L_2;H_{02})$. In either way, one eventually ends up with an $A_\infty$-category, which we denote by $\F(M)$.

\subsection{Characteristic cycles}
Given a real-analytic manifold $Z$, one can consider sheaves of
$\K$-vector spaces which are constructible (with respect to some real
analytic stratification, which may depend on the sheaf). Denote by
$D_c(Z)$ the full subcategory of
the bounded derived
category of sheaves of $\K$-vector spaces comprising complexes with constructible cohomology. Kashiwara's characteristic cycle construction \cite{Kashiwara} associates to any object ${\mathcal G}$ in this category a Lagrangian cycle $CC({\mathcal G})$ inside $M = T^*Z$, which is a cone (invariant under rescaling of cotangent fibres). If ${\mathcal G}$ is the structure sheaf of a closed submanifold, this cycle is just the conormal bundle, but otherwise it tends to be singular. Nadler and Zaslow consider the structure sheaves of submanifolds $W \subset Z$ which are (real-analytic but) not necessarily closed. For each such ``standard object'', they construct a smoothing of $CC(\mathcal G)$, which is a Lagrangian submanifold of $M$. In the special case where $W$ is an open subset with smooth boundary, this is essentially equivalent to the construction indicated in Example \ref{Ex:cotangent}. The singular boundary case is considerably more complicated, and leads to Lagrangian submanifolds which are generally only asymptotically invariant under the Liouville flow (their limits at infinity are singular Legendrian cycles). Still, one can use them as objects of a Fukaya-type category, which is a variant of the previously described construction. We denote it by $\A(M)$, where $\A$ stands for ``asymptotic''. The main result of \cite{NZ} is

\begin{Theorem} \label{th:nz1}
The smoothed characteristic cycle construction gives rise to a full embedding of derived categories, $D_c(Z) \longrightarrow D\A(M)$.
\end{Theorem}

The proof relies on two ideas. One of them, namely, that the standard objects generate $D_c(Z)$, is more elementary (it can be viewed as a fact about decompositions of real subanalytic sets). For purely algebraic reasons, this means that it is enough to define the embedding only on standard objects. The other, more geometric, technique is the reduction of Floer cohomology to Morse theory, provided by the work of Fukaya-Oh \cite{FOh}.

\subsection{Decomposing the diagonal\label{subsec:nadler-equi}}

In $D_c(Z)$, the structure sheaves of any two distinct points are
algebraically disjoint (there are no morphisms between them). Of
course, the same holds for cotangent fibres in the Fukaya category of
$M$, which are the images of such structure sheaves under the
embedding from Theorem \ref{th:nz1}. As a consequence, in $\F(M)$ (or
$\A(M)$, the difference being irrelevant at this level) one cannot
expect to have a finite resolution of a closed exact $L \subset M$ in
terms of cotangent fibres. However, using the Wehrheim-Woodward
formalism of Lagrangian correspondences \cite{Wehrheim-Woodward},
Nadler proves a modified version of this statement, where the fibres are replaced by standard objects associated to certain contractible subsets of $Z$.

Concretely, fix a real analytic triangulation of $Z$. Denote by $x_i$ the
vertices of the triangulation, by $U_i$ their stars, and by $U_I =
\bigcap_{i \in I} U_i$ the intersections of such stars, indexed by finite
sets $I = \{i_0,\dots,i_d\}$. There is a standard Cech resolution of the
constant sheaf $\K_Z$ in terms of the $\K_{U_I}$ (we will encounter a
similar construction again later on, in Section \ref{subsec:cech}). Rather
than applying this to $Z$ itself, we take the diagonal inclusion $\delta:
Z \rightarrow Z \times Z$, and consider the induced resolution of
$\delta_*(\K_Z)$ by the objects $\delta_*(\K_{U_I})$. Consider the
embedding from Theorem \ref{th:nz1} applied to $Z \times Z$. The image of
$\delta_*(\K_Z)$ is the conormal bundle of the diagonal $\Delta =
\delta(Z)$, and each $\delta_*(\K_{U_I})$ maps to the standard object
associated to $\delta(U_I) \subset Z \times Z$. Since each $U_I$ is
contractible, one can deform $\delta(U_I)$ to $U_I \times \{x_{i_d}\}$ (as
locally closed submanifolds of $Z \times Z$), and this induces an isotopy
of the associated smoothed characteristic cycles. \emph{A priori}, this isotopy
is not compactly supported, hence not well-behaved in our category (it
does not preserve the isomorphism type of objects). However, this is not a
problem if one is only interested in morphisms from or to a given closed
Lagrangian submanifold.

To formalize this, take $\A^{cpt}(M)$ to be the subcategory of $\A(M)$
consisting of closed Lagrangian submanifolds (this is in fact the Fukaya
category in the classical sense). Dually, let $mod(\A(M))$,
$mod(\A^{cpt}(M))$ be the associated categories of $A_\infty$-modules.
There is a chain of $A_\infty$-functors
\begin{equation} \label{eq:yoneda}
\A^{cpt}(M) \longrightarrow \A(M) \longrightarrow mod(\A(M))
\longrightarrow mod(\A^{cpt}(M)).
\end{equation}
The first and second one, which are inclusion and the Yoneda embedding,
are full and faithful. The last one, restriction of $A_\infty$-modules,
will not generally have that property. However, the composition of all
three is just the Yoneda embedding for $\A^{cpt}$, which is again full and
faithful. In view of \cite{Wehrheim-Woodward}, and its chain-level
analogue \cite{MWW}, each object $C$ in $\A(M
\times M)$ (and more generally, twisted complex built out of such objects)
induces a convolution functor
\begin{equation}
\Phi_C : \A^{cpt}(M) \longrightarrow mod(\A^{cpt}(M))
\end{equation}
(usually, one reverses the sign of the symplectic form on one of the two
factors in $M \times M$, but for cotangent bundles, this can be
compensated by a fibrewise reflection $\sigma: M \rightarrow M$). First,
take $C$ to be the conormal bundle of the diagonal, which is the same as
the graph of $\sigma$. Then, convolution with $C$ is isomorphic to the
embedding \eqref{eq:yoneda}. On the other hand, if $C$ is the smoothed
characteristic cycle of some product $U \times \{x\}$, then $\Phi_C$ maps
each object to a direct sum of copies of $T^*_x$ (the image of
that fibre under the functor $\A(M) \rightarrow mod(\A^{cpt}(M))$, to
be precise). Finally, if $C$ is just a Lagrangian submanifold, $\Phi_C$ is
invariant under Lagrangian isotopies which are not necessarily compactly
supported. By combining those facts, one obtains the desired resolution of
a closed exact $L \subset M$. Nadler actually pushes these ideas somewhat
further, using a refined version of this argument, to show that:

\begin{Theorem} \label{th:nz2}
The embedding $D_c(Z) \rightarrow D\A(M)$ from Theorem \ref{th:nz1} is an
equivalence.
\end{Theorem}

\begin{Remark} \label{Rmk:cechdiff}
If one is only interested in the spectral sequence \eqref{eq:ss}, there
may be a potential simplification, which would bypass some of the
categorical constructions above. First of all, rewrite $HF^*(L,L) \iso
H^*(L;\K)$ as $HF^*(L \times L,\nu^*\Delta)$. Then, using the resolution
of $\nu^*\Delta$ in $\A(M \times M)$ described above, one gets a spectral
sequence converging towards that group, whose $E_1$ page comprises the Floer
cohomology groups between $L \times L$ and the smoothed characteristic
cycles of $\delta(U_I)$. Since $L \times L$ is compact, one can deform
$\delta(U_I)$ to $U_I \times \{x_{i_d}\}$, and then further isotop its
smoothed characteristic cycle to $T^*_{x_{i_d}} \otimes T^*_{x_{i_d}}$. In
the terminology used in \eqref{eq:ss}, this brings the terms in the $E_1$
page into the form $End(E_L)_{x_{i_d}}$. To get the desired $E_2$ term,
one would further have to check that the differentials reproduce the ones
in the Cech complex with twisted coefficients in $End(E_L)$. This of
course follows from Theorem \ref{th:nz2}, but there ought to be a more
direct geometric argument, just by looking at the relevant spaces of
holomorphic triangles; this seems a worth while endeavour, but we have not
attempted to study it in detail.
\end{Remark}

\section{Lefschetz thimbles\label{sec:vc}}

This section gives an overview of the paper \cite{FSS}, and an account -- emphasising geometric rather than algebraic aspects -- of some of the underlying theory from the book \cite{FCPLT}.

\subsection{Fukaya categories of Lefschetz fibrations}

In principle, the notion of Lefschetz fibration can be defined in a purely symplectic way. However, we will limit ourselves to the more traditional algebro-geometric context. Let $X$ be a smooth affine variety, and
\begin{equation}
\pi: X \longrightarrow \C
\end{equation}
a polynomial, which has only nondegenerate critical points. For convenience, we assume that no two such points lie in the same fibre. Additionally, we impose a condition which excludes singularities at infinity, namely: let $\bar{X}$ be a projective completion of $X$, such that $D = \bar{X} \setminus X$ is a divisor with normal crossing. We then require that (for an approriate choice of $\bar{X}$) the closure of $\pi^{-1}(0)$ should be smooth in a neighbourhood of $D$, and intersect each component of $D$ transversally. Finally, for Floer-theoretic reasons, we require $X$ to be Calabi-Yau (have trivial canonical bundle).

Take any K{\"a}hler form on $\bar{X}$ which comes from a metric on ${\mathcal O}(D)$. Its restriction to $X$ makes that variety into a Weinstein manifold (of finite type, but not complete; the latter deficiency can, however, be cured easily, by attaching the missing part of the conical end). Moreover, parallel transport for $\pi$ is well-defined away from the singular fibres, in spite of its non-properness.

A vanishing path $\gamma:[0,\infty) \rightarrow \C$  is an embedding starting at a critical value $\gamma(0)$ of $\pi$, and such that for $t \gg 0$, $\gamma(t) = {\text{\it const.}} - it$ is a half-line going to $-i\infty$. To each such path one can associate a \emph{Lefschetz thimble} $\Delta_\gamma \subset X$, which is a Lagrangian submanifold diffeomorphic to $\R^n$, projecting properly to $\gamma([0,\infty)) \subset \C$. More precisely, $\gamma^{-1} \circ \pi|\Delta_{\gamma}$ is the standard proper Morse function on $\R^n$ with a single minimum (placed at the unique critical point of $\pi$ in the fibre over $\gamma(0)$). When defining the Floer cohomology between two Lefschetz thimbles, the convention is to rotate the semi-infinite part of the first path in anticlockwise direction for some small angle. Omitting certain technical points, this can be interpreted as adding a Hamiltonian term as in Section \ref{subsec:epsilon}. In particular, one again has
\begin{equation} \label{eq:vc-1}
HF^*(L_\gamma,L_\gamma) \iso H^*(L_\gamma;\K) = \K.
\end{equation}
Now suppose that $(\gamma_0,\dots,\gamma_m)$ is a basis (sometimes also called a distinguished basis) of vanishing paths. We will not recall the definition here; for a sketch, see Figure \ref{fig:basis}. In that situation, if one takes $(\gamma_j,\gamma_k)$ with $j>k$ and applies the rotation described above to $\gamma_j$, the result remains disjoint from $\gamma_k$. Hence,
\begin{equation} \label{eq:vc-2}
HF^*(L_{\gamma_j},L_{\gamma_k}) = 0 \;\; \text{for all $j>k$.}
\end{equation}
\begin{figure}[ht]
\begin{centering}
\setlength{\unitlength}{0.00066667in}
\begingroup\makeatletter\ifx\SetFigFont\undefined%
\gdef\SetFigFont#1#2#3#4#5{%
  \reset@font\fontsize{#1}{#2pt}%
  \fontfamily{#3}\fontseries{#4}\fontshape{#5}%
  \selectfont}%
\fi\endgroup%
{\renewcommand{\dashlinestretch}{30}
\begin{picture}(2862,1870)(0,-10)
\put(1000,312){\makebox(0,0)[lb]{{\SetFigFont{10}{12}{\rmdefault}{\mddefault}{\updefault}$\gamma_1$}}}
\put(2250,1212){\circle*{70}}
\put(1275,1812){\circle*{70}}
\put(1650,1212){\circle*{70}}
\dashline{60.000}(2325,312)(2325,12)
\dashline{60.000}(2850,312)(2850,12)
\dashline{60.000}(1275,312)(1275,12)
\drawline(1275,612)(1275,312)
\dashline{60.000}(225,312)(225,12)
\drawline(1050,1212)(1053,1211)(1058,1209)
	(1068,1206)(1084,1200)(1106,1193)
	(1133,1184)(1166,1173)(1204,1160)
	(1245,1145)(1289,1130)(1334,1115)
	(1380,1099)(1425,1083)(1469,1068)
	(1511,1053)(1550,1039)(1587,1026)
	(1622,1014)(1655,1002)(1685,990)
	(1713,980)(1740,970)(1764,960)
	(1787,951)(1809,942)(1830,933)
	(1850,924)(1877,913)(1903,901)
	(1928,889)(1952,876)(1975,864)
	(1998,851)(2020,839)(2041,826)
	(2062,813)(2081,799)(2099,786)
	(2116,773)(2132,760)(2147,748)
	(2160,735)(2172,723)(2184,710)
	(2194,698)(2204,686)(2213,674)
	(2223,659)(2233,644)(2242,627)
	(2251,611)(2259,593)(2267,576)
	(2274,558)(2281,540)(2287,523)
	(2292,506)(2297,490)(2301,475)
	(2304,461)(2308,448)(2310,436)
	(2313,424)(2316,408)(2318,394)
	(2320,380)(2322,366)(2323,351)
	(2324,337)(2324,325)(2325,316)
	(2325,313)(2325,312)
\drawline(2250,1212)(2252,1210)(2257,1206)
	(2266,1198)(2280,1186)(2298,1170)
	(2321,1150)(2347,1128)(2375,1103)
	(2405,1076)(2436,1049)(2466,1022)
	(2495,996)(2522,972)(2548,949)
	(2571,927)(2593,907)(2613,889)
	(2631,871)(2647,855)(2662,840)
	(2675,826)(2688,813)(2700,799)
	(2716,780)(2732,762)(2746,744)
	(2759,726)(2771,707)(2782,690)
	(2792,672)(2800,655)(2808,638)
	(2815,622)(2821,606)(2825,591)
	(2829,577)(2832,563)(2835,550)
	(2838,537)(2840,522)(2842,507)
	(2843,491)(2845,474)(2846,454)
	(2847,433)(2848,410)(2849,386)
	(2849,363)(2850,342)(2850,326)
	(2850,317)(2850,313)(2850,312)
\drawline(1275,1812)(1273,1811)(1267,1809)
	(1257,1806)(1241,1800)(1220,1793)
	(1193,1783)(1160,1771)(1123,1758)
	(1082,1743)(1038,1727)(993,1710)
	(949,1694)(904,1677)(861,1660)
	(821,1644)(782,1629)(746,1614)
	(712,1600)(681,1586)(652,1573)
	(625,1560)(600,1548)(577,1535)
	(556,1523)(536,1511)(518,1499)
	(500,1487)(479,1471)(459,1455)
	(440,1438)(422,1421)(405,1402)
	(388,1383)(373,1363)(358,1343)
	(345,1321)(332,1299)(321,1277)
	(310,1253)(301,1230)(292,1206)
	(284,1182)(277,1158)(271,1134)
	(266,1110)(261,1085)(257,1061)
	(253,1037)(250,1012)(247,991)
	(245,969)(243,946)(241,922)
	(239,897)(238,871)(236,843)
	(235,813)(233,780)(232,746)
	(231,709)(230,671)(229,630)
	(228,589)(228,547)(227,506)
	(227,466)(226,429)(226,397)
	(225,369)(225,347)(225,331)
	(225,320)(225,315)(225,312)
\drawline(1275,612)(1275,613)(1275,617)
	(1275,626)(1274,641)(1274,659)
	(1273,677)(1272,695)(1271,711)
	(1269,726)(1267,739)(1265,750)
	(1263,762)(1259,774)(1256,785)
	(1251,798)(1245,810)(1238,824)
	(1231,837)(1222,850)(1212,864)
	(1201,876)(1189,889)(1176,900)
	(1163,912)(1151,921)(1139,929)
	(1126,938)(1112,948)(1097,957)
	(1081,967)(1065,977)(1048,987)
	(1031,997)(1014,1007)(998,1017)
	(982,1026)(967,1036)(952,1045)
	(938,1053)(925,1062)(911,1072)
	(897,1082)(883,1092)(870,1103)
	(857,1115)(845,1127)(834,1139)
	(823,1151)(814,1164)(806,1176)
	(800,1188)(795,1200)(790,1212)
	(788,1224)(786,1235)(785,1247)
	(785,1258)(785,1271)(787,1283)
	(790,1296)(793,1310)(798,1323)
	(803,1336)(810,1349)(817,1361)
	(825,1372)(833,1383)(842,1394)
	(852,1403)(863,1412)(874,1420)
	(886,1428)(899,1436)(913,1443)
	(929,1450)(945,1457)(962,1463)
	(980,1469)(999,1475)(1017,1480)
	(1036,1484)(1054,1488)(1073,1492)
	(1090,1495)(1108,1497)(1125,1499)
	(1142,1502)(1160,1503)(1178,1505)
	(1196,1506)(1215,1507)(1233,1508)
	(1253,1508)(1272,1509)(1290,1508)
	(1308,1508)(1326,1507)(1342,1506)
	(1358,1505)(1373,1503)(1387,1502)
	(1400,1499)(1417,1496)(1433,1493)
	(1449,1488)(1465,1483)(1480,1476)
	(1495,1469)(1509,1461)(1522,1453)
	(1533,1443)(1544,1433)(1554,1423)
	(1563,1412)(1570,1402)(1577,1391)
	(1584,1378)(1591,1363)(1598,1346)
	(1607,1327)(1615,1306)(1624,1283)
	(1632,1261)(1640,1241)(1645,1226)
	(1648,1217)(1650,1213)(1650,1212)
\put(0,537){\makebox(0,0)[lb]{{\SetFigFont{9}{10.8}{\rmdefault}{\mddefault}{\updefault}$\gamma_0$}}}
\put(2025,312){\makebox(0,0)[lb]{{\SetFigFont{9}{10.8}{\rmdefault}{\mddefault}{\updefault}$\gamma_2$}}}
\put(2900,537){\makebox(0,0)[lb]{{\SetFigFont{9}{10.8}{\rmdefault}{\mddefault}{\updefault}$\gamma_3$}}}
\put(1050,1212){\circle*{70}}
\end{picture}
}
\caption{\label{fig:basis}}
\end{centering}
\end{figure}

As usual, rather than working on the level of Floer cohomology, we
want to have underlying $A_\infty$-structures. There is a convenient
shortcut, which eliminates noncompact Lagrangian submanifolds from the
foundations of the theory (but which unfortunately requires $char(\K)
\neq 2$). Namely, let $\tilde{X} \rightarrow X$ be the double cover
branched over some fibre $\pi^{-1}(-iC)$, $C \gg 0$. Roughly speaking,
one takes the ordinary Fukaya category ${\mathcal F}(\tilde{X})$,
which contains only compact Lagrangian submanifolds, and defines
$\F(\pi)$ to be its $\Z/2$-invariant part (only invariant objects and
morphisms; obviously, getting this to work on the cochain level
requires a little care). This time, let's allow only vanishing paths
which satisfy $\gamma(t) = -it$ for $t \geq C$ (which is no problem,
since each path can be brought into this form by an isotopy). One then
truncates the Lefschetz thimble associated to such a path, so that it
becomes a Lagrangian disc with boundary in $\pi^{-1}(-iC)$, and takes
its preimage in $\tilde{X}$, which is a closed $\Z/2$-invariant
Lagrangian sphere $\tilde{S}_\gamma \subset \tilde{X}$. On the
cohomological level, the $\Z/2$-invariant parts of the Floer
cohomologies of these spheres still satisfy the same properties as
before, in particular reproduce (\ref{eq:vc-1}) and (\ref{eq:vc-2}).

\begin{Theorem}
If $(\gamma_0,\dots,\gamma_m)$ is a basis of vanishing paths, the associated $\tilde{S}_{\gamma_j}$ form a full exceptional collection in the derived category $D\F(\pi)$.
\end{Theorem}

The fact that we get an exceptional collection is elementary; it just
reflects the two equations \eqref{eq:vc-1} and \eqref{eq:vc-2} (or
rather, their counterparts for the modified definition of Floer
cohomology involving double covers). Fullness, which is the property
that this collection generates the derived category, is rather more
interesting. The proof given in \cite{FCPLT} relies on the fact that
the product of Dehn twists along the $\tilde{L}_{\gamma_k}$ is
isotopic to the covering involution in $\tilde{X}$. Hence, if $L \subset X$ is a closed Lagrangian submanifold which lies in $\pi^{-1}(\{\im(z) > -C\})$, this product of Dehn twists will exchange the two components of the preimage $\tilde{L} \subset \tilde{X}$. The rest of the argument essentially consists in applying the long exact sequence from \cite{Se}.

\subsection{Postnikov decompositions}
We will use some purely algebraic properties of exceptional
collections, see for instance \cite{Goro} (the subject has a long
history in algebraic geometry; readers interested in this might find
the collection \cite{Rudakov} to be a good starting point). Namely,
let $C$ be a triangulated category, linear over a field $\K$, and let
$(Y_0,\dots,Y_m)$ be a full exceptional collection of objects in $C$. Then, for any object $X$, there is a collection of exact triangles
\begin{equation} \label{eq:tower}
Z_k \otimes Y_k \rightarrow X_k \rightarrow X_{k-1}
\stackrel{[1]}{\longrightarrow} Z_k\otimes Y_k
\end{equation}
where $X_m = X$, $X_{-1} = 0$, and $Z_k = Hom_C^*(Y_k, X_k)$ (morphisms of all degrees; by assumption, this is finite-dimensional, so $Z_k \otimes Y_k$ is the direct sum of finitely many shifted copies of $Y_k$). The map $Z_k \otimes Y_k \rightarrow X_k$ is the canonical evaluation map, and $X_{k-1}$ is defined (by descending induction on $k$) to be its mapping cone. To get another description of $Z_k$, one can use the unique (right) Koszul dual exceptional collection $(Y_m^!,\dots,Y_0^!)$, which satisfies
\begin{equation} \label{eq:ss2}
Hom_C^*(Y_j,Y_k^!) = \begin{cases} \K \text{ (concentrated in degree zero)} & j = k, \\ 0 & \text{otherwise}. \end{cases}
\end{equation}
It then follows by repeatedly applying \eqref{eq:tower} that $Z_k^\vee \iso Hom_C^*(X,Y_k^!)$. Now, given any cohomological functor $R$ on $C$, we get an induced spectral sequence converging to $R(X)$, whose starting page has columns $E_1^{rs} = (Z^{\vee}_{m-r} \otimes R(Y_{m-r}))^{r+s}$. In particular, taking $R = Hom_C^*(-,X)$ and using the expression for $Z_k$ explained above, we get a spectral sequence converging to $Hom_C^*(X,X)$, which starts with
\begin{equation}
E_1^{rs} = \big(Hom_C^*(X,Y_{m-r}^!) \otimes Hom_C^*(Y_{m-r},X)\big)^{r+s}.
\end{equation}

We now return to the concrete setting where $C = D\F(\pi)$. In this case, the Koszul dual of an exceptional collection given by a basis of Lefschetz thimbles is another such basis $\{\gamma_0^!,\dots,\gamma_m^!\}$.
This is a consequence of the more general relation between mutations (algebra) and Hurwitz moves on vanishing paths (geometry). Applying \eqref{eq:ss2}, and going back to the original definition of Floer cohomology, we therefore get the following result: for every (exact, graded, spin) closed Lagrangian submanifold $L \subset M$, there is a spectral sequence converging to $HF^*(L,L) \iso H^*(L;\K)$, which starts with
\begin{equation} \label{eqn:specseq}
E_1^{rs} = (HF(L, \Delta_{\gamma_{m-r}^!}) \otimes HF(\Delta_{\gamma_{m-r}},L))^{r+s}.
\end{equation}

\subsection{Real algebraic approximation\label{subsec:nt}}
The existence of \eqref{eqn:specseq} is a general statement about Lefschetz fibrations. To make the connection with cotangent bundles, we use a form of the Nash-Tognoli theorem, see for instance \cite{NT}, namely:

\begin{Lemma} \label{Lem:13}
If $Z$ is a closed manifold and $p :Z \rightarrow \R$ is a Morse function, there is a Lefschetz fibration $\pi: X \rightarrow \C$ with a compatible real structure, and a diffeomorphism $f: Z \rightarrow X_{\R}$, such that $\pi \circ f$ is $C^2$-close to $p$.
\end{Lemma}

The diffeomorphism $f$ can be extended to a symplectic embedding
$\phi$ of a neighbourhood of the zero-section of $Z \subset M = T^*Z$
into $X$. Hence (perhaps after a preliminary radial rescaling) we can transport closed exact Lagrangian submanifolds $L \subset M$ over to $X$. The critical points of $\pi$ fall into two classes, namely real and purely complex ones, and the ones in the first class correspond canonically to critical points of $p$. By a suitable choice of vanishing paths, one can ensure that
\begin{equation} \label{eqn:infibre}
\begin{aligned}
 & HF^*(\Delta_{\gamma_k},\phi(L)) \\ & \iso \begin{cases} HF^*(T^*_{x_k},L) &
 \text{if $\gamma(0) \in X_\R$ corresponds to $x_k \in Crit(p)$}, \\
 0 & \text{otherwise}.
\end{cases}
\end{aligned}
\end{equation} 
Here, the Floer cohomology on the right hand side is taken inside
$T^*Z$. The same statement holds for the other groups in
\eqref{eqn:specseq}, up to a shift in the grading which depends on the
Morse index of $x_k$. As a consequence, the starting page of that
spectral sequence can be thought of as $C^*_{Morse}(Z;End^*(E_L))$,
where the Morse complex is taken with respect to the function $p$, and
using the (graded) local coefficient system $E_L$. This is one page
earlier than our usual starting term \eqref{eq:ss}, but is already
good enough to derive Theorem \ref{Thm:main} by appealing to some classical
manifold topology (after taking the product with a sphere if
necessary, one can assume that $dim(Z)>5$, in which case simple
connectivity of $Z$ implies that one can choose a Morse function
without critical points of index or co-index $1$). 

\begin{Remark}
The differential on the $E_1$-page of  (\ref{eqn:specseq}) is given in \cite[Corollary
18.27]{FCPLT}, in terms of holomorphic triangle products between
adjacent Lefschetz thimbles in the exceptional collection.  In the
special situation of (\ref{eqn:infibre}), identifying the $E_1$ page
with $C^*_{Morse}(Z;End(E_L))$, there is also the Morse differential
$\delta$ coming
from parallel transport in the local system $E_L \rightarrow Z$
(compare Remark \ref{Rmk:cechdiff}).  For Lefschetz fibrations arising from
real algebraic approximation, rather than some more canonical
construction, there seems to be no reason for these to agree in
general; but one does expect the parts of the 
differential leading out of the first column, and into the last column,
to agree.  For instance, by deforming the final vanishing path to lie along
the real axis, one can ensure that the entire thimble leading out of
the maximum $x_{max}$ of the Morse function is contained in the real
locus, after which the intersection points between this thimble and
one coming from a critical point $x$ of index one less correspond
bijectively to the gradient lines of the Morse function between
$x_{max}$ and $x$ (the situation at the minimum $x_{min}$ is
analogous).  If it was known that this part of the
$E_1$-differential did reproduce the corresponding piece of $\delta$,
one could hope to study non-simply-connected cotangent bundles in this
approach. 

\end{Remark}

\section{Family Floer cohomology\label{sec:family}}
This section covers the third point of view on Theorem
\ref{Thm:main}. This time, the presentation is less linear, and
occasionally several ways of reaching a particular goal are
sketched. The reader should keep in mind the preliminary nature of
this discussion. In some parts, this means that there are complete but
unpublished constructions. For others, only outlines or strategies of
proof exist, in which case we will be careful to formulate the relevant statements as conjectures.

\subsection{Cech complexes\label{subsec:cech}}

At the start of the paper, we mentioned that to an $L \subset M =
T^*Z$ one can associate the bundle $E_L$ of Floer cohomologies $E_x =
HF^*(T^*_x,L)$. One naturally wants to replace $E_L$ by an underlying
cochain level object ${\mathcal E}_L$, which should be a ``sheaf of
complexes'' in a suitable sense. In our interpretation, this will be a
dg module over a dg algebra of Cech cochains (there are several other
possibilities, with varying degrees of technical difficulty; see
\cite[Section 5]{Fukaya} and \cite[Section 4]{Nadler} for sketches of
two of these).

Fix a smooth triangulation of $Z$, with vertices $x_i$, and let $U_I$ be the intersections of sets in the associated open cover, just as in Section \ref{subsec:nadler-equi} (but omitting the real analyticity condition). This time, we want to write down the associated Cech complex explicitly, hence fix an ordering of the $i$'s. Let $\Gamma(U_I)$ be the space of locally constant $\K$-valued functions on $U_I$, which in our case is $\K$ if $U_I \neq \emptyset$, and $0$ otherwise. The Cech complex is
\begin{equation}
{\mathcal C} = \bigoplus_I \Gamma(U_I)[-d]
\end{equation}
where $d = |I|-1$. This carries the usual differential, and also a
natural associative product making it into a dg algebra. Namely, for
every possible splitting of $I'' = \{i_0 < \dots < i_d\}$ into $I' =
\{i_0 < \cdots < i_k\}$ and $I = \{i_k < \cdots < i_d\}$, one takes
\begin{equation} \label{eq:cech-product}
\Gamma(U_I) \otimes \Gamma(U_{I'}) \xrightarrow{\text{restriction}}
\Gamma(U_{I''}) \otimes \Gamma(U_{I''}) \xrightarrow{\text{multiplication}} \Gamma(U_{I''}).
\end{equation}
We want to consider (unital right) dg modules over ${\mathcal C}$. Denote the dg category of such modules by ${\mathcal M} = mod({\mathcal C})$. This is \emph{not} a dg model for the derived category: there are acyclic modules which are nontrivial in $H({\mathcal M})$, and as a consequence, quasi-isomorphism does not imply isomorphism in that category.

All objects we will consider actually belong to a more restricted class, distinguished by a suitable ``locality'' property; we call these dg modules of {\em presheaf type}. The definition is that such a dg module ${\mathcal E}$ needs to admit a splitting
\begin{equation} \label{eq:presheaf}
{\mathcal E} = \bigoplus_I {\mathcal E}_I[-d],
\end{equation}
where the sum is over all $I = \{i_0 < \cdots < i_d\}$ such that $U_I
\neq \emptyset$. This splitting is required to be compatible with the
differential and module structure. This means that the differential
maps ${\mathcal E}_I$ to the direct sum of ${\mathcal E}_{I'}$ over
all $I' \supset I$; that $1 \in \Gamma(U_{i'})$ acts as the identity
on ${\mathcal E}_I[-d]$ for all $I = \{i_0 < \cdots < i_d\}$ with $i_d
= i'$, and as zero otherwise; and that the component ${\mathcal E}_I
\otimes \Gamma(U_{I'}) \longrightarrow {\mathcal E}_{I''}$ of the
module structure can only be nonzero if $I'' \supset I \cup I'$. In
particular, the ``stalks'' ${\mathcal E}_I$ themselves are
subquotients of ${\mathcal E}$, and inherit a dg module structure from
that.   The stalks associated to the smallest subsets of $Z$ (i.e. to
maximal index sets $I$) are actually chain complexes, with all chain
homomorphisms being module endomorphisms, from which one can show:

\begin{Lemma} \label{Th:presheaf}
Let ${\mathcal E}$ be a dg module of presheaf type. If each ${\mathcal
  E}_I$ is acyclic, ${\mathcal E}$ itself is isomorphic to the zero
object in $H({\mathcal M})$.
\end{Lemma}

Here is a first example. Let $P \rightarrow Z$ be a flat $\K$-vector
bundle (or local coefficient system). For each $I = \{i_0 < \cdots <
i_d\}$ such that $U_I \neq \emptyset$, define $({\mathcal E}_P)_I =
P_{x_{i_d}}$. Form the direct sum ${\mathcal E}_P$ as in
\eqref{eq:presheaf}. Equivalently, one can think of this as the sum of
$P_{x_{i_d}} \otimes \Gamma(U_I)$ over all $I$, including empty
ones. The differential consists of terms $({\mathcal E}_P)_I
\rightarrow ({\mathcal E}_P)_{I'}$ for $I = \{i_0 < \cdots < i_d\}$, $I'
= I \cup \{i'\}$. If $i'<i_d$, these are given (up to sign) by
restriction maps $P_{x_{i_d}} \otimes \Gamma(U_I) \rightarrow
P_{x_{i_d}} \otimes \Gamma(U_{I'})$. In the remaining nontrivial case
where $i'>i_d$ and $U_{I'} \neq \emptyset$, there is a unique edge of
our triangulation going from $x_{i_d}$ to $x_{i'}$, and one combines
restriction with parallel transport along that edge. The module
structure on ${\mathcal E}_P$ is defined in the obvious way, following
the model \eqref{eq:cech-product};  compatibility of the differential
and the right $\mathcal{C}$-module structure is ensured by our choice $({\mathcal E}_P)_I =
P_{x_{i_d}}$, taking the fibre of $P$ over the vertex corresponding to
the last index $i_d$ of $I$. Clearly, $H^*({\mathcal E}_P) \iso
H^*(Z;P)$ is ordinary cohomology with $P$-coefficients. Moreover,
using the fact that every ${\mathcal E}_P$ is free as a module
(ignoring the differential), one sees that
\begin{equation} \label{eq:pp}
H^*(hom_{\mathcal M}({\mathcal E}_{P_0},{\mathcal E}_{P_1})) \iso H^*(Z;Hom(P_0,P_1)).
\end{equation}
In particular, referring back to the remark made above, this leads to quite concrete examples of acyclic dg modules which are nevertheless nontrivial objects.

Still within classical topology, but getting somewhat closer to the intended construction, suppose now that $Q \rightarrow Z$ is a differentiable fibre bundle, equipped with a connection. For $I = {i_0< \cdots < i_d}$ with $U_I \neq \emptyset$, define
\begin{equation} \label{eq:cubical}
({\mathcal E}_Q)_I = C_{-*}(Q_{x_{i_d}}),
\end{equation}
where $C_{-*}$ stands for cubical cochains, with the grading reversed (to go with our general cohomological convention). As before, we want to turn the (shifted) direct sum of these into a dg module over ${\mathcal C}$. The module structure is straightforward, but the differential is a little more interesting, being the sum of three terms. The first of these is the ordinary boundary operator on each $({\mathcal E}_Q)_I$. The second one is the Cech differential $({\mathcal E}_Q)_I \rightarrow ({\mathcal E}_Q)_{I'}$, where $I = \{i_0 < \cdots < i_d\}$ and $I' = I \cup \{i'\}$ with $i' < i_d$. The final term consists of maps
\begin{equation} \label{eq:family-cont}
C_{-*}(Q_{x_{i_d}}) \longrightarrow C_{-*}(Q_{x_{i'_e}})[1-e]
\end{equation}
where $I = \{i_0 < \cdots < i_d\}$, $I' = I \cup \{i'_1 < \cdots < i'_e\}$ with $i_d < i'_1$. Take the standard $e$-dimensional simplex $\Delta^e$. It is a classical observation \cite{Adams} that there is a natural family of piecewise smooth paths in $\Delta^e$, parametrized by an $(e-1)$-dimensional cube, which join the first to the last vertex. In our situation, the triangulation of $Z$ contains a unique simplex with vertices $\{i_d,i'_1,\dots,i'_e\}$, and we therefore get a family of paths joining $x_{i_d}$ to $x_{i'_e}$. The resulting parametrized parallel transport map
\begin{equation}
[0;1]^{e-1} \times Q_{x_{i_d}} \longrightarrow Q_{x_{i'_e}}
\end{equation}
induces \eqref{eq:family-cont} (up to sign; the appearance of a cube here motivated our use of cubical chains). It is not difficult to show that the resulting total differential on ${\mathcal E}_Q$ indeed squares to zero, and is compatible with the module structure.

Finally, let's turn to the symplectic analogue of this, in which one starts with a closed exact Lagrangian submanifold $L \subset M = T^*Z$ (subject to the usual conditions: if one wants $\Z$-graded modules, $L$ should be graded; and if one wants to use coefficient fields $char(\K) \neq 2$, it should be relatively spin). The appropriate version of \eqref{eq:cubical} is
\begin{equation}
({\mathcal E}_L)_I = CF^*(T^*_{x_{i_d}},L)
\end{equation}
and again, ${\mathcal E}_L$ is the sum of these. In the differential,
we now use the Floer differential on each $CF^*$ summand (replacing
the differential on cubical chains), and continuation maps or their
parametrized analogues (rather than parallel transport maps), which
govern moving one cotangent fibre into the other. (Related continuation
maps appeared in \cite{Kasturirangan-Oh2}.)  Of course, the details are somewhat different from the previous case. To get a chain map $CF^*(T^*_{x_{i_d}},L) \rightarrow CF^*(T^*_{x_{i'_1}},L)$, one needs a path
\begin{equation}
\gamma: \R \longrightarrow Z, \quad \gamma(s) = x_{i_d}\text{ for $s \ll 0$, } \gamma(s) = x_{i'_1} \text{ for $s \gg 0$.}
\end{equation}
More generally, families of such paths parametrized by $[0;1)^{e-1}$ appear. In the limit as one (or more) of the parameters go to $1$, the path breaks up into two (or more) pieces separated by increasingly long constant stretches. This ensures that the usual composition laws for continuation maps apply, compare \cite{Salamon-Zehnder}. Still, with these technical modifications taken into account, the argument remains essentially the same as before.

\subsection{Wrapped Floer cohomology\label{subsec:wrapped}}
One naturally wants to extend the previous construction to allow $L$ to be noncompact (for instance, a cotangent fibre). This would be impossible using the version of Floer cohomology from Section \ref{subsec:epsilon}, since that does not have sufficiently strong isotopy invariance properties: $HF^*(T^*_x,L)$ generally depends strongly on $x$. Instead, we use a modified version called ``wrapped Floer cohomology''. This is not fundamentally new: it appears in \cite{AS} for the case of cotangent bundles, and is actually the open string analogue of Viterbo's symplectic cohomology \cite{Vit}.

Take a Weinstein manifold $M$ (complete and of finite type), and
consider exact Lagrangian submanifolds inside it which are Legendrian
at infinity; this time, no real analyticity conditions will be
necessary. We again use Hamiltonian functions of the form $H(r,y) =
h(e^r)$ at infinity, but now require that $\lim_{r \rightarrow \infty}
h'(e^r) = +\infty$. This means that as one goes to infinity, the
associated Hamiltonian flow is an unboundedly accelerating version of
the Reeb flow. Denote by $CW^*(L_0,L_1) = CF^*(L_0,L_1;H) =
CF^*(\phi_X^1(L_0),L_1)$ the resulting Floer complex, and by
\begin{equation}
HW^*(L_0,L_1) = H(CW^*(L_0,L_1))
\end{equation}
its cohomology, which we call {\em wrapped Floer cohomology}. This is
independent of the choice of $H$. Moreover, it remains invariant under
isotopies of either $L_0$ or $L_1$ inside the relevant class of
submanifolds. Such isotopies need  no longer be compactly supported; the Legendrian submanifolds at infinity may move (which is exactly the property we wanted to have). By exploiting this, it is easy to define a triangle product on wrapped Floer cohomology. In the case where $L_0 = L_1 = L$, we have a natural map from the ordinary cohomology $H^*(L)$ to $HW^*(L,L)$, which however is generally neither injective nor surjective. For instance, for $L = \R^n$ inside $M = \C^n$, $HW^*(L,L)$ vanishes. Another, far less trivial, example is the following one:

\begin{Theorem}[Abbondandolo-Schwarz] \label{Pro:2}
Let $M=T^*Z$ be the cotangent bundle of a closed oriented manifold, and take two cotangent fibres
$L_0 = T_{x_0}^*Z$, $L_1 =T_{x_1}^*Z$.  Then $HW^*(L_0,L_1) \iso H_{-*}(\mathcal{P}_{x_0,x_1})$ is the (negatively graded) homology of the space of paths in $Z$ going from $x_0$ to $x_1$.
\end{Theorem}

This is proved in \cite{AS}; the followup paper \cite{AS2} shows that this isomorphism sends the triangle product on $HW^*$ to the Pontryagin product (induced by composition of paths) on path space homology.

We want our wrapped Floer cochain groups to carry an $A_\infty$-structure, refining the cohomology level product. When defining this, one encounters the same difficulty as in \eqref{eq:strict-product}. Again, there is a solution based on a rescaling trick: one takes $h(t) = \frac{1}{2} t^2$, and uses the fact that $\phi_X^2$ differs from $\phi_Y^{-\log(2)} \phi_X^1 \phi_Y^{\log(2)}$ ($\phi_Y$ being the Liouville flow) only by a compactly supported isotopy. This is particularly intuitive for cotangent bundles, where the radial coordinate at infinity is $e^r = |p|$; then, $H = h(e^r) = \half|p|^2$ gives rise to the standard geodesic flow $(\phi_X^t)$, and $Y = p\partial_p$ is the rescaling vector field, meaning that $\phi_Y^{\log(2)}$ doubles the length of cotangent vectors. An alternative (not identical, but ultimately quasi-isomorphic) approach to the same problem is to define wrapped Floer cohomology as a direct limit over functions $H_k$ with more moderate growth $H_k(r,y) = ke^r$ at infinity. However, the details of this are quite intricate, and we will not describe them here. In either way, one gets an $A_\infty$-category $\W(M)$, which we call the {\em wrapped Fukaya category}. For cotangent bundles $M = T^*Z$, a plausible cochain level refinement of Theorem \ref{Pro:2} is the following

\begin{Conjecture} \label{Conj:abbondandolo-schwarz}
Let $L = T^*_x$ be a cotangent fibre. Then, the $A_\infty$-structure on $CW^*(L,L)$ should be quasi-isomorphic to the dg algebra structure on $C_{-*}(\Omega_x)$, where $\Omega_x$ is the (Moore) based loop space of $(Z,x)$.
\end{Conjecture}

Returning to the general case, we recall that a fundamental property of symplectic cohomology, established in \cite{Vit}, is {\em Viterbo functoriality} with respect to embeddings of Weinstein manifolds. One naturally expects a corresponding property to hold for wrapped Fukaya categories. Namely, take a bounded open subset $U \subset M$ with smooth boundary, such that $Y$ points outwards along $\partial U$. One can then attach an infinite cone to $\partial U$, to form another Weinstein manifold $M' = U \cup_{\partial U} ([0;\infty) \times \partial U)$. Suppose that $(L_1,\dots,L_d)$ is a finite family of exact Lagrangian submanifolds in $M$, which are Legendrian at infinity, and with the following additional property: $\theta|L_k = dR_k$ for some compactly supported function $R_k$, which in addition vanishes in a neighbourhood of $L_k \cap \partial U$. This implies that $L_k \cap \partial U$ is a Legendrian submanifold of $\partial U$. Again, by attaching infinite cones to $L_k \cap U$, one gets exact Lagrangian submanifolds $L_k' \subset M'$, which are Legendrian at infinity. Let $\A \subset \W(M)$ be the full $A_\infty$-subcategory with objects $L_k$, and similarly $\A' \subset \W(M')$ for $L_k'$. Then,

\begin{Conjecture}[Abouzaid-Seidel] \label{Conj:viterbo}
There is a natural (up to isomorphism) $A_\infty$-functor ${\mathcal R}: \A \rightarrow \A'$.
\end{Conjecture}

Note that, even though we have not mentioned this explicitly so far, all $A_\infty$-categories under consideration have units (on the cohomological level), and ${\mathcal R}$ is supposed to be a (cohomologically) unital functor. Hence, the conjecture implies that various relations between objects, such as isomorphism or exact triangles, pass from $\A$ to $\A'$, which is a nontrivial statement in itself.

\subsection{Family Floer cohomology revisited\label{subsec:rev}} 
For $M = T^*Z$, there is a straightforward variation of the
construction from Section \ref{subsec:cech}, using wrapped Floer
cohomology instead of ordinary Floer cohomology. This associates to any exact Lagrangian submanifold $L \subset M$, which is Legendrian at infinity, a dg module ${\mathcal E}_L$ over ${\mathcal C}$. In fact, it gives rise to an $A_\infty$-functor
\begin{equation} \label{eq:g-functor}
{\mathcal G}: \W(M) \longrightarrow {\mathcal M} = mod({\mathcal C}).
\end{equation}
While little has been rigorously proved so far, there are plausible expectations for how this functor should behave, which we will now formulate precisely. Take $Q \rightarrow Z$ to be the path fibration (whose total space is contractible). Even though this is not strictly a fibre bundle, the construction from Section \ref{subsec:cech} applies, and yields a dg module ${\mathcal E}_Q$, which has $H^*({\mathcal E}_Q) \iso \K$ (this can be viewed as a resolution of the simple ${\mathcal C}$-module).

\begin{Conjecture} \label{Conj:faith}
Let $L = T^*_x$ be a cotangent fibre. Then, ${\mathcal E}_L$ is isomorphic to ${\mathcal E}_Q$ in $H({\mathcal M})$. Moreover, if $Z$ is simply-connected, ${\mathcal G}$ gives rise to a quasi-isomorphism $C_{-*}(\Omega_x) \iso CW^*(L,L) \rightarrow hom_{\mathcal M}({\mathcal E}_L,{\mathcal E}_L)$.
\end{Conjecture}

The first statement can be seen as a parametrized extension of Theorem
\ref{Pro:2}. A possible proof would be to consist of checking that the
chain level maps constructed in \cite{AS} can be made compatible with
parallel transport (respectively continuation) maps, up to a suitable
hierarchy of chain homotopies. This would yield a map of dg modules;
to prove that it is an isomorphism, one would then apply Lemma
\ref{Th:presheaf} to its mapping cone. The second part of the
conjecture is less intuitive. The assumption of simple connectivity is necessary, since otherwise the endomorphisms of ${\mathcal E}_Q$ may not reproduce the homology of the based loop space (see Section \ref{subsec:algebra} for further discussion of this); however, it is not entirely clear how that would enter into a proof.

\begin{Conjecture} \label{Conj:generate}
Any one fibre $L = T^*_x$ generates the derived category (taken, as usual, to be the homotopy category of twisted complexes) of $\W(M)$.
\end{Conjecture}

There are two apparently quite viable approaches to this, arising from the contexts of Sections \ref{sec:nz} and \ref{sec:vc}, respectively. To explain the first one, we go back to the general situation where $M$ is a finite type complete Weinstein manifold, whose end is modelled on the contact manifold $N$, and where real analyticity conditions are imposed on $N$ and its Legendrian submanifolds. Then, if $(L_0,L_1)$ are exact Lagrangian submanifolds which are Legendrian at infinity, we have a natural homomorphism
\begin{equation}
HF^*(L_0,L_1) \longrightarrow HW^*(L_0,L_1),
\end{equation}
which generalizes the map $H^*(L) \rightarrow HW^*(L,L)$ mentioned in
Section \ref{subsec:wrapped}. These maps are compatible with triangle
products, and even though the details have not been checked, it seems
plausible that they can be lifted to an $A_\infty$-functor $\F(M)
\rightarrow \W(M)$. Actually, what one would like to use is a variant
of this, where $\F(M)$ is replaced by the Nadler-Zaslow category
$\A(M)$, or at least a sufficiently large full subcategory of
it. Assuming that this can be done, one can take the generators of
$\A(M)$ provided by the proof of Theorem \ref{th:nz2}, and then map
them to $\W(M)$, where isotopy invariance ought to ensure that they
all become isomorphic to cotangent fibres (note that in the wrapped Fukaya category, any two cotangent fibres are mutually isomorphic).

The second strategy for attacking Conjecture \ref{Conj:generate} is fundamentally similar, but based on Lefschetz fibrations. Recall that with our definition, the total space of a Lefschetz fibration $\pi: X \rightarrow \C$ is itself a finite type Weinstein manifold. One then expects to have a canonical functor $\F(\pi) \rightarrow \W(X)$. Given a Lefschetz fibration constructed as in Section \ref{subsec:nt} by complexifying a Morse function on $Z$, one would then combine $\F(\pi) \rightarrow \W(X)$ with the functor from Conjecture \ref{Conj:viterbo}, applied to a small neighbourhood of $Z$ embedded into $X$. The outcome would be that the restrictions of Lefschetz thimbles generate the wrapped Fukaya category of that neighbourhood. One can easily check that all such restrictions are isomorphic to cotangent fibres.

Suppose now that $Z$ is simply-connected. In that case, if one accepts Conjectures \ref{Conj:faith} and \ref{Conj:generate}, it follows by purely algebraic means that \eqref{eq:g-functor} is full and faithful. Take the dg module ${\mathcal E}_L = \bigoplus_I ({\mathcal E}_L)_I[-d]$ obtained from a closed exact $L \subset M$, and equip it with the decreasing filtration by values of $d = |I|-1$. The associated graded space is precisely the dg module ${\mathcal E}_P$ constructed from the local coefficient system $P_x = HW^*(T^*_x,L) = HF^*(T^*_x,L)$. Hence, one gets a spectral sequence which starts with $H^*(Z;End(P))$ according to \eqref{eq:pp}, and converges to the group $H(hom_{\mathcal M}({\mathcal E}_L,{\mathcal E}_L))$, which would be equal to $HF^*(L,L) \iso H^*(L;\K)$ as a consequence of the previously made conjecture. This explains how \eqref{eq:ss} arises from this particular approach.


\subsection{The (co)bar construction\label{subsec:algebra}}

There is a more algebraic perspective, which provides a shortcut through part of the argument above. To explain this, it is helpful to recall the classical bar construction. Let ${\mathcal C}$ be a dg algebra over our coefficient field $\K$, with an augmentation $\varepsilon: {\mathcal C} \rightarrow \K$, whose kernel we denote by ${\mathcal I}$. One can then equip the free tensor coalgebra $T({\mathcal I}[1])$ with a differential, and then dualize it to a dg algebra ${\mathcal B} = T({\mathcal I}[1])^\vee$. Consider the standard resolution ${\mathcal R} = T({\mathcal I}[1]) \otimes {\mathcal C}$ of the simple ${\mathcal C}$-module $\mathcal C/\mathcal I$. One can prove that the endomorphism dga of ${\mathcal R}$ as an object of $mod(\mathcal{C})$ is quasi-isomorphic to ${\mathcal B}$. For standard algebraic reasons, this induces a quasi-equivalence between the subcategory of $mod(\mathcal{C})$ generated by ${\mathcal R}$, and the triangulated subcategory of $mod(\mathcal B)$ generated by the free module $\mathcal B$. Denote that category by $modf(\mathcal B)$.

The relevance of this duality to our discussion is a basic connection to loop spaces, which goes back to Adams \cite{Adams}. He observed that if $Z$ is simply-connected, and $\mathcal C$ is the dg algebra of Cech cochains, then ${\mathcal B}$ is quasi-isomorphic to $C_{-*}(\Omega Z)$. Hence, reversing the functor above, we get a full and faithful embedding
\begin{equation}
H(modf(C_{-*}(\Omega Z))) \longrightarrow H(\mathcal M),
\end{equation}
where $\mathcal M = mod(\mathcal C)$ as in \eqref{eq:g-functor}. If
one moreover  assumes that Conjectures \ref{Conj:abbondandolo-schwarz} and \ref{Conj:generate} hold, it follows that $\W(M)$ itself is derived equivalent to $H(modf(C_{-*}(\Omega Z)))$. Hence, one would get a full embedding of the wrapped Fukaya category into $H(\mathcal M)$ for algebraic reasons, avoiding the use of Cech complexes altogether.

\section{The non-simply-connected case}

In this final section, we discuss exact Lagrangian submanifolds in non-simply-connected cotangent bundles. Specifically, we prove Corollary \ref{Cor:main}, and then make a few more observations about the wrapped Fukaya category.

\subsection{A finite covering trick\label{subsec:trick}}

We start by recalling the setup from Section \ref{sec:nz}. Fix a real analytic structure on $Z$, and the associated category ${\mathcal A}(M)$. For any closed exact Lagrangian submanifold $L\subset M=T^*Z$ which is spin and has zero Maslov class, we have the spectral sequence
\[
E_2^{pq} = H^p(Z;End^q(E_L)) \Rightarrow H^*(L;\K)
\]
arising from the resolution of $L$ by ``standard objects'' in
${\mathcal A}(M)$. From now on, we assume that $char(\K) = p>0$. In
that case, one can certainly find a finite covering $b: \tilde{Z}
\rightarrow Z$ such that $b^*E_L$ is trivial (as a graded bundle of
$\K$-vector spaces). In fact, $E_L$ gives rise to a representation
$\rho: \pi_1(Z) \rightarrow GL_r(\K)$, and one takes $\tilde{Z}$ to be
the Galois covering associated to $ker(\rho) \subset \pi_1(Z)$. Set
$\tilde{M} = T^*\tilde{Z}$, denote by $\beta: \tilde{M} \rightarrow M$
the induced covering, and by $\tilde{L} = \beta^{-1}(L)$ the preimage
of our Lagrangian submanifold. This inherits all the properties of
$L$, hence we have an analogous spectral sequence for
$\tilde{L}$. Note also that for  obvious reasons, the associated bundle of Floer cohomologies satisfies
\begin{equation} \label{eq:Lem1}
E_{\tilde{L}} \iso b^*E_L,
\end{equation}
hence is trivial by definition. Now, the discussion after
\eqref{eq:ss} carries over with no problems to the non-simply-connected
case if one assumes that the Lagrangiaan submanifold is connected, and
its bundle of Floer cohomologies is trivial. In particular, one gets
that $\tilde{L}$ is in fact connected: if it weren't, the analogue of
Theorem \ref{Thm:main}(iii) would apply to its connected components
(since their Floer cohomology bundles are also trivial), implying that
any two would have to intersect each other, which is a
contradiction. With connectedness at hand, it then follows by the same
argument that \eqref{eq:Lem1} has one-dimensional fibres. Hence,
$End(E_L)$ is trivial, which means that the spectral sequence for the
original $L$ degenerates, yielding $H^*(L;\K) \iso H^*(Z;\K)$. In
fact, by borrowing arguments from \cite{Nadler} or from \cite{FSS},
one sees that in the Fukaya category of $M$, $L$ is isomorphic to the
zero-section equipped with a suitable $spin$ structure (the difference
between that and the \emph{a priori} chosen spin structure on $Z$ is
precisely described by the bundle $E_L$, which, note, has structure
group $\pm 1$). Using that, one also gets the analogue of part (iii)
of Theorem \ref{Thm:main}.  Finally, the cohomological restrictions
(i)-(iii) of Theorem \ref{Thm:main} for arbitrary fields of positive
characteristic imply them 
for $\K$ with $char(\K)=0$, which completes the proof of Corollary
\ref{Cor:main}. 


\subsection{The Eilenberg-MacLane case}

For general algebraic reasons, there is an $A_\infty$-functor
\begin{equation} \label{eqn:algembed}
\W(M) \longrightarrow mod(CW^*(T^*_x,T^*_x)).
\end{equation}
Here, the right hand side is the dg category of $A_\infty$-modules over
the endomorphism $A_\infty$-algebra of the object $T^*_x$. Now suppose
that $Z = K(\Gamma,1)$ is an Eilenberg-MacLane space, so that the
cohomology of $CW^*(T^*_x,T^*_x)$ is concentrated in degree zero. By the
Homological Perturbation Lemma, this implies that the $A_\infty$-structure
is formal, hence by Theorem \ref{Pro:2} quasi-isomorphic to the group
algebra $\K\,[\Gamma]$ (this argument allows us to avoid Conjecture
\ref{Conj:abbondandolo-schwarz}). However, in order to make
\eqref{eqn:algembed} useful, we do want to appeal to
(the currently unproven) Conjecture \ref{Conj:generate}.
Assuming from now on that that is true, one finds that
\eqref{eqn:algembed} is a full embedding, and actually lands (up to
functor isomorphism) in the subcategory $modf(\K\,[\Gamma])$ generated by
the one-dimensional free module. As far as that subcategory is concerned,
one could actually replace $A_\infty$-modules by ordinary chain complexes
of $\K\,[\Gamma]$-modules. This gives a picture of $\W(M)$ in classical
algebraic terms.

To see how this might be useful, let's drop the assumption that the Maslov
class is zero, and consider general closed, exact, and
spin Lagrangian submanifolds $L \subset M$. Neither of the two arguments
in favour of Conjecture \ref{Conj:generate} sketched in Section
\ref{subsec:rev} actually uses the Maslov index. It is therefore plausible
to assume that the description of $\W(M)$ explained above still
applies. Denote by 
${\mathcal N}_L$ the $\Z/2$-graded $A_\infty$-module over $\K\,[\Gamma]$
corresponding to $L$, and by $N_L$ its cohomology module. In fact, $N_L$
is nothing other than our previous family Floer cohomology bundle $E_L$,
now considered as a module over $\Gamma = \pi_1(Z)$. There is a purely
algebraic obstruction theory, which determines when an $A_\infty$-module
is formal (isomorphic to its cohomology). The obstructions lie in
\begin{equation} \label{eq:ext}
Ext^r_{\K[\Gamma]}(N_L,N_L[1-r]), \;\; r \geq 2
\end{equation}
where the $[1-r]$ now has to be interpreted mod $2$. In particular, if all
those groups vanish, it would follow directly that ${\mathcal N}_L$ is
formal. However, there is no particular \emph{a priori} reason why that should
happen.

Returning to the trick from Section \ref{subsec:trick}, assume now that
$char(\K) = p>0$. Then, after passing to a finite cover, one can assume
that $N_L$ is a direct sum of trivial representations, hence
\eqref{eq:ext} is a direct sum of copies of $H^r(\Gamma;\K)$. One can try
to kill the relevant obstructions by passing to further finite covers.
Generally, this is unlikely to be successful (there are examples of mod
$p$ cohomology classes which survive pullback to any finite cover, see for
instance \cite[Theorem 6.1]{FW}). However, in the special case where $H^*(\Gamma;\K)$
is generated by degree $1$ classes, such as for surfaces or tori, it is
obviously possible. In those cases, one could then find a finite cover $b:
\tilde{Z} \rightarrow Z$, inducing $\beta:\tilde{M}\rightarrow M$, such that the $A_\infty$-module ${\mathcal
N}_{\tilde{L}}$ associated to $\tilde{L} = \beta^{-1}(L)$ is isomorphic to a
direct sum of ordinary trivial modules. Just by looking at the
endomorphism ring of this object, it becomes clear that there can actually
be only one summand, so $\tilde{L}$ is isomorphic to the zero-section. One
can then return to $M$ by the same argument as before, and obtain the same
consequences as in the original Theorem \ref{Thm:main}. This would be a
potential  application
of the machinery from Section \ref{sec:family} which has no obvious
counterpart in the other approaches.  

Finally, note that the appeal to Conjecture \ref{Conj:generate} above can be sidestepped, at least working over fields of characteristic not equal to two, in the special case where $T^*Z$ admits a Lefschetz pencil for which, for suitable vanishing paths, the Lefschetz thimbles are \emph{globally} cotangent fibres.  This is, of course, an exceptional situation, but it can be achieved for suitable Lefschetz pencils on the affine variety $T^*T^n = (\C^*)^n$ (complexify a Morse function which is a sum of height functions on the distinct $S^1$-factors of $T^n$).

\bibliographystyle{alpha}

\end{document}